\documentclass[12pt]{article}
\usepackage{amsmath,amsfonts,amssymb,amsthm,amscd}

\input xy
\xyoption{all}

\title{Model theory and K\"ahler geometry}
\author{Rahim Moosa\thanks{Supported by an NSERC grant}\\University of Waterloo \and Anand Pillay\thanks{Supported by a Marie Curie chair}\\University of  Leeds}
\date{17 November, 2006}

\newtheorem{Theorem}{Theorem}[section]
\newtheorem{Proposition}[Theorem]{Proposition}
\newtheorem{Definition}[Theorem]{Definition} 
\newtheorem{Remark}[Theorem]{Remark}

\newtheorem{Corollary}[Theorem]{Corollary}
\newtheorem{Fact}[Theorem]{Fact}

\theoremstyle{definition}
\newtheorem{Example}[Theorem]{Example}

\newcommand{\R}{\mathbb R}   
  
\newcommand{\Z}{\mathbb Z}

\newcommand{\C}{\mathbb C}
\newcommand{\Pp}{\mathbb P}

\newcommand{\A}{\mathcal A}

\newcommand{\re}{\operatorname{Re}}
\newcommand{\im}{\operatorname{Im}}
\newcommand{\vol}{\operatorname{vol}}
\newcommand{\Hom}{\operatorname{Hom}}
\newcommand{\Th}{\operatorname{Th}}
\newcommand{\acl}{\operatorname{acl}}

\begin{document}
\maketitle

\begin{abstract} We survey and explain some recent work at the intersection of model theory and bimeromorphic geometry (classification of compact complex manifolds). Included here are the essential saturation of the many-sorted structure ${\cal C}$ of K\"ahler manifolds, the conjectural role of hyperk\"ahler manifolds in the description of strongly minimal sets in ${\cal C}$, and Campana's work on the isotriviality of hyperk\"ahler families and its connection with the nonmultidimensionality conjecture.

\end{abstract}
\section{Introduction}

The aim of this paper is to discuss in some detail the relationship between ideas from model theory (classification theory, geometric stability theory) and  those from bimeromorphic geometry (classification of compact complex manifolds), with reference to current research.
Earlier work along these lines is in \cite{pillay2000}, \cite{pillay01}, \cite{pillayscanlon2001}, \cite{pillayscanlon2000}, \cite{radin2004}, \cite{ret}, \cite{moosa-ccs}, \cite{sat}, \cite{ams}.
We will also take the opportunity here to describe for model-theorists some of the basic tools of complex differential geometry, as well as summarise important notions, facts and theorems such as the Hodge decomposition, and local Torelli.

 Zilber \cite{zilber93} observed some time ago that if a compact complex manifold $M$ is considered naturally as a first order structure (with predicates for analytic subsets of $M$, $M\times M$, etc.) then $\operatorname{Th}(M)$ has finite Morley rank. The same holds if we consider the category ${\cal A}$ of compact complex (possibly singular) spaces as a many-sorted first order structure. This observation of Zilber was closely related, historically, to the work on Zariski structures and geometries by Hrushovski and Zilber \cite{hz93}. 
 
There is a rich general theory of theories of finite Morley rank, encompassing both Shelah's work on classification theory (classifying first order theories and their models) as well as the more self-conciously geometric theory of $1$-basedness (modularity), definable groups, definable automorphism groups, etc...  It turns out that $\operatorname{Th}({\cal A})$ witnesses most of the richness of this theory.  Among the main points of the current article is that notions belonging to the Shelah theory such as nonorthogonality and nonmultidimensionality, have a very clear geometric content, and are connected with things such as ``variation of Hodge structure". 
 
The class of compact K\"ahler manifolds has been identified as an important rather well-behaved class of compact complex manifolds, where there is a better chance of classification.  The first author~\cite{sat} observed that such manifolds can to all intents and purposes be treated as saturated structures (inside which one can apply the compactness theorem).  We give some more details in section~\ref{saturationsection}, explaining the role of the K\"ahler condition.

The category of compact K\"ahler manifolds (or rather compact complex analytic spaces that are holomorphic images of compact K\"ahler manifolds) is a ``full reduct" of the many-sorted structure ${\cal A}$. We call it ${\cal C}$. In \cite{pillayscanlon2001} it was pointed out how, from work of Lieberman, one can see that $\operatorname{Th}({\cal A})$ is about as complicated as it can be from the point of view of Shelah's theory (it has the DOP). We have conjectured on the other hand that $\operatorname{Th}({\cal C})$ is rather tame. $\operatorname{Th}({\cal C})$ could not be uncountably categorical (unidimensional),
 but we believe it to be the next best thing, nonmultidimensional. The description of $U$-rank $1$ types (equivalently {\em simple} compact complex manifolds) in $\operatorname{Th}({\cal C})$ which are {\em trivial} is still open, and it is conjectured that they are closely related to so-called irreducible hyperk\"ahler manifolds. As we explain in section 5, an isotriviality result for families of hyperk\"ahler manifolds in ${\cal C}$, due to Campana, represents some confirmation of the nonmultdimensionality of $\operatorname{Th}({\cal C})$. 
 
 We now give a brief survey of the model theory of compact complex manifolds, continuing in a sense \cite{pillay-thisvolume}.
 There are several published survey-type articles, such as \cite{moosa-ccs} and \cite{pillay2000}, to which the interested reader is referred for more details.
 We assume familiarity with the notion of a complex manifold $M$.
 An {\em analytic} subset $X$ of $M$ is a subset $X$ such that for each $a\in M$ there is an an open neighbourhood $U$ of $a$ in $M$ such that $X\cap U$ is the common zero-set of a finite set of holomorphic functions on $U$.  A compact complex manifold $M$ is viewed as a first order structure by equipping it with predicates for all the analytic subsets of $M$ and its cartesian powers.
 The fundamental fact observed by Zilber is that the theory of this first order structure has quantifier elimination and has finite Morley rank.  We can of course consider the collection of {\em all}  compact complex manifolds (up to biholomorphism) and view it is a many-sorted first order structure (predicates for all analytic subsets of cartesian products of sorts) which again has QE and finite Morley rank (sort by sort). The same holds for the larger class of {\em compact complex analytic spaces}.  A compact complex analytic space is (a compact topological space) locally modelled on zero-sets of finitely many holomorphic functions on open domains in $\C^{n}$ with of course biholomorphic transition maps. We have the notion of an analytic subset of $X$ and its cartesian powers, and we obtain thereby a first order structure as before. As above we let ${\cal A}$ denote the many-sorted structure of compact complex analytic spaces, and we let $\mathcal{L}$ denote its language.  (A complex analytic space is usually presented as a ringed space, where the rings may be nonreduced. We refer to \cite{moosa-ccs} for more discussion of this. In any case by a {\em complex variety}  we will mean a reduced and irreducible complex analytic space.)
 
If $X$ is a compact complex variety and $a\in X$, then $\{a\}$ is an analytic subset of $X$, and hence is essentially named by a constant. So really ${\cal A}$ has names for all elements (of all sorts).  Let ${\cal A}'$ be a saturated elementary extension of ${\cal A}$. If $S$ is a sort in ${\cal A}$, we let $S'$ be the corresponding sort in ${\cal A}'$. For example, $(\Pp_{1})'$ denotes the projective line over a suitable elementary extension $\C '$ of $\C$.

Among the basic facts connecting definability and geometry are: 

(i) For $X,Y$ sorts in ${\cal A}$ the definable maps from $X$ to $Y$ are precisely the piecewise meromorphic maps.

(ii) If $p(x)$ is a complete type of $\operatorname{Th}({\cal A})$ then $p$ is the generic type (over ${\cal A}$) of a unique compact complex variety $X$.
That is, $p(x)$ is axiomatised by ``$x\in X'$ but $x\notin Y'$ for any proper analytic  subset $Y\subset X$".

(iii) If $a, b$ are tuples from ${\cal A}'$ with $tp(a)$ the generic type of $X$ and $tp(b)$ the generic type of $Y$, then $dcl(a) = dcl(b)$ iff $X$ and $Y$ are bimeromorphic. 

(iv) Suppose $a,b$ are tuples from ${\cal A}'$, and $tp(a/b)$ is stationary. Then there are compact complex varieties $X,Y$ and a meromorphic dominant map $f:X\to Y$  whose fibres over a non-empty Zariski open subset of $Y$ are irreducible, and such that: $ab$ is a generic point (realizes the generic type) of $X$, $b$ is a generic point of $Y$, and (in ${\cal A}'$) $f(ab) = b$.
So $tp(a/b)$ is the ``generic type" of the ``generic fibre" $X_{b}$ of $f:X\to Y$.
We consider definable sets such as $X_{b}$ as ``nonstandard"  analytic subsets of $X$.

Algebraic geometry lives in ${\cal A}$ on the sort $\Pp_{1}$.
Any irreducible complex quasi-affine {\em algebraic} variety $V$ has a compactification ${\bar V}$ which will be a compact complex variety living as a sort in ${\cal A}$, and biholomorphic with a closed subvariety of $\Pp_1^n$ for some $n>0$.
The variety $V$ will be a Zariski open, hence definable, subset of ${\bar V}$.

A compact complex variety $X$ is said to be Moishezon if $X$ is bimeromorphic with a complex projective algebraic variety.
This is equivalent to $X$ being internal to the sort $\Pp_{1}$, and also equivalent to a generic point $a$ of $X$ being in the definable closure of elements from $(\Pp_{1})'$.  The expression ``algebraic" is sometimes used in place of Moishezon.  We extend naturally this notion to nonstandard analytic sets as well as definable sets and stationary types in ${\cal A}'$. The ``strong conjecture" from \cite{pillay-thisvolume} then holds in ${\cal A}'$ in the more explicit form: if ($Y_{z}:z\in Z)$ is a normalized family of definable subsets of a definable set $X$, then for $a\in X$, $Z_{a} = \{z\in Z: a$ is generic on $Y_{z}\}$ is Moishezon  (namely generically internal to $(\Pp_{1})'$).  This result was derived in \cite{pillay01} from a theorem due independently to Campana and Fujiki. 

In \cite{pillayscanlon2001} it was shown that any  strongly minimal modular group definable in ${\cal A}$ is definably isomorphic to a complex torus. An appropriate generalization to strongly minimal modular groups in ${\cal A}'$ was obtained in \cite{ams}. 

More details on the classification of strongly minimal sets (or more generally types of $U$-rank $1$) in $\Th({\cal C})$ will appear in section 5.

\bigskip

\section{Preliminaries on complex forms}
\label{preliminaries}

We give in this section a brief review of the basic theory of complex-valued differential forms.
The reader may consult~\cite{fritzsche-grauert02} or~\cite{voisin} for a more detailed treatment of this material.

Suppose $X$ is an $n$-dimensional complex manifold.
By a {\em co-ordinate system} $(z,U)$ on $X$ we mean an open set $U\subset X$ and a homeomorphism $z$ from $U$ to an open ball in $\C^n$.
Composing with the co-ordinate projections we obtain {\em complex co-ordinates} $z_i:U\to \C$ for $i=1,\dots,n$, which we decompose into real and imaginary parts as $z_i=x_i+iy_i$.

Fix a co-ordinate system $(z,U)$ on $X$ and a point $x\in U$.
Let $T_{X,x}$ denote the (real) tangent space of $X$ at $x$.
Viewed as the space of $\R$-linear derivations on real-valued smooth functions at $x$, we have that
$$\left\{\frac{\partial}{\partial x_1},\dots,\frac{\partial}{\partial x_n},\frac{\partial}{\partial y_1},\dots,\frac{\partial}{\partial y_n}\right\}$$
forms an $\R$-basis for $T_{X,x}$.
But the complex manifold structure on $X$ gives $T_{X,x}$ also an $n$-dimensional {\em complex} vector space structure, which can be described as follows.
Let $T_{X,\C,x}:=T_{X,x}\otimes\C$ denote the complexification of the real tangent space.
We have a  decomposition $T_{X,\C,x}=T_{X,x}^{1,0}\oplus T_{X,x}^{0,1}$ where $T_{X,x}^{1,0}$ is the complex subspace generated by
$$\left\{\frac{\partial}{\partial z_i}:=\frac{\partial}{\partial x_i}-i\frac{\partial}{\partial y_i} \mid i=1,\dots,n\right\}$$
and $T_{X,x}^{0,1}$ is generated by
$$\left\{\frac{\partial}{\partial \overline{z_i}}:=\frac{\partial}{\partial x_i}+i\frac{\partial}{\partial y_i} \mid i=1,\dots,n\right\}.$$
If we view $T_{X,\C,x}$ as the space of $\C$-linear derivations on {\em complex}-valued smooth functions at $x$, then $T_{X,x}^{1,0}$ corresponds to those that vanish on all the anti-holomorphic functions (functions whose complex conjugates are holomorphic at $x$).
In any case, $T_{X,x}^{1,0}$ is called the {\em holomorphic tangent space} of $X$ at $x$.
The natural inclusion $T_{X,x}\subset T_{X,\C,x}$ followed by the projection $T_{X,\C,x}\to T_{X,x}^{1,0}$ produces an $\R$-linear isomorphism between the real tangent space and the holomorphic tangent space.
This isomorphism makes $T_{X,x}$ canonically into a complex vector space.

Despite our presentation, the above constructions do not depend on the co-ordinates and extend globally:
We have the complexification of the (real) tangent bundle $T_{X,\C}:= T_X\otimes\C$, and a decomposition into complex vector subbundles, $T_{X,\C}=T_X^{1,0}\oplus T_X^{0,1}$, whereby the holomorphic tangent bundle   $T_X^{1,0}$ is naturally isomorphic as a real vector bundle with $T_X$.
It is with respect to this isomorphism that we treat $T_X$ as a complex vector bundle.

A {\em complex-valued differential $k$-form} (or just a {\em $k$-form}) at $x\in X$ is an alternating $k$-ary $\R$-multilinear map $\phi:T_{X,x}\times\cdots\times T_{X,x}\to \C$.
The complex vector space of all $k$-forms at $x$ is denoted by $F^k_{X,\C,x}$.
The {\em real} differential $k$-forms, $F^k_{X,\R,x}$, are exactly those forms in $F_{X,\C,x}^k$ that are real-valued.
So $F_{X,\C,x}^k=F^k_{X,\R,x}\otimes\C$.
In particular,
$F_{X,\C,x}^1=\operatorname{Hom}_{\R}(T_{X,x},\R)\otimes\C$ is the complexification of the real cotangent space at $x$.
Hence, in a co-ordinate system $(z,U)$ about $x$, if we let $\left\{dx_1,\dots,dx_n,dy_1,\dots,dy_n\right\}$ be the dual basis to $\left\{\frac{\partial}{\partial x_1},\dots,\frac{\partial}{\partial x_n},\frac{\partial}{\partial y_1},\dots,\frac{\partial}{\partial y_n}\right\}$ for $\operatorname{Hom}_{\R}(T_{X,x},\R)$, then
$$dz_i:=dx_i+idy_i,$$
$$d\overline{z}_i:=dx_i-idy_i$$
for $i=1,\dots,n$, form a $\C$-basis for $F_{X,\C,x}^1$.

Now $F_{X,\C,x}^k$ is the $k$th exterior power of $F_{X,\C,x}^1$.
Given $I=(i_1,\dots,i_p)$ an increasing sequence of numbers between $1$ and $n$, let $dz_I$ be the $p$-form $dz_{i_1}\wedge\cdots\wedge dz_{i_p}$.
Similarly, let $d\overline{z_I}:=d\overline{z_{i_1}}\wedge\cdots\wedge d\overline{z_{i_p}}$.
Then 
$$\left\{dz_I\wedge d\overline{z_J}\mid I=(i_1,\dots,i_p), J=(j_1,\dots,j_q), p+q=k\right\}$$
is a $\C$-basis for $F_{X,\C,x}^k$.
This gives us a natural decomposition
$$F_{X,\C,x}^k=\bigoplus_{p+q=k}F_{X,x}^{p,q}$$
where $F_{X,x}^{p,q}$ is generated by the forms $dz_I\wedge d{\overline{z_J}}$ where $I=(i_1,\dots,i_p)$ and $J=(j_1,\dots,j_q)$ are increasing sequences of numbers between $1$ and $n$.
The complex vector subspaces $F_{X,x}^{p,q}$ can also be more intrinsically described as made up of those $k$-forms $\phi$ such that
$$\phi(cv_1,\dots,cv_n)=c^p\overline{c}^q\phi(v_1,\dots,v_n)$$
for all $v_1,\dots,v_k\in T_{X,x}$ and $c\in\C$.
Such forms are said to be of {\em type} $(p,q)$.

Once again, these constructions extend globally to $X$ and we have complex vector bundles $\displaystyle F_{X,\C}^k=\bigoplus_{p+q=k}F_{X,\C}^{p,q}$.
For $U\subseteq X$ an open set, by a complex $k$ form {\em on $U$} we mean a smooth section to the bundle $F_{X,\C}^k$ over the set $U$.
Similarly for forms of {\em type $(p,q)$ on $U$}.
We denote by $\mathcal{A}^k$ and $\mathcal{A}^{p,q}$ the sheaves on $X$ of smooth sections to $F_{X,\C}^k$  and $F_{X,\C}^{p,q}$ respectively.
So $\mathcal{A}^k(U)$ is the space of all complex $k$-forms on $U$ while $\mathcal{A}^{p,q}(U)$ is the space of all complex forms of type $(p,q)$ on $U$.
Given a co-ordinate system $(z,U)$ on $X$, a $k$-form $\omega\in\mathcal{A}^k(U)$ can be expressed as
$\displaystyle \omega=\sum_{|I|+|J|=k}f_{IJ}dz_I\wedge d\overline{z_J}$
where $f_{IJ}:U\to\C$ are smooth.
Note that $dz_i$ and $d\overline{z_i}$ are being viewed here as $1$-forms on $U$.
By convention, $\mathcal{A}^0$ is the sheaf of $\C$-valued smooth functions.

The
{\em exterior derivative} map $d:\mathcal{A}^k(U)\to\mathcal{A}^{k+1}(U)$ is defined by
$$d\left(\sum_{|I|+|J|=k}f_{IJ}dz_I\wedge d\overline{z_J}\right)
= \sum_{|I|+|J|=k}df_{IJ}\wedge dz_I\wedge d\overline{z_J}$$
where for any smooth function $f:U\to\C$, $df\in \mathcal{A}^1(U)$ is given by
$$df:=\sum_{i=1}^n\frac{\partial f}{\partial z_i}dz_i+\sum_{i=1}^n\frac{\partial f}{\partial \overline{z_i}}d\overline{z_i}.$$
If we define
$\displaystyle \partial f:=\sum_{i=1}^n\frac{\partial f}{\partial z_i}dz_i$
and
$\displaystyle \overline{\partial}f=\sum_{i=1}^n\frac{\partial f}{\partial \overline{z_i}}d\overline{z_i},$
and then extend these maps so that
$\partial:\mathcal{A}^{p,q}(U)\to\mathcal{A}^{p+1,q}(U)$ is given by 
$$\partial\left(\sum f_{IJ}dz_I\wedge d\overline{z_J}\right)
= \sum\partial f_{IJ}\wedge dz_I\wedge d\overline{z_J}$$
and $\overline{\partial}:\mathcal{A}^{p,q}(U)\to\mathcal{A}^{p,q+1}(U)$ is given by 
$$\overline{\partial}\left(\sum f_{IJ}dz_I\wedge d\overline{z_J}\right)
= \sum\overline{\partial} f_{IJ}\wedge dz_I\wedge d\overline{z_J};$$
then we see that $d=\partial+\overline{\partial}$.

One can show that $d$, $\partial$, and $\overline{\partial}$ are all independent of the co-ordinate system and extend to sheaf maps on $\mathcal{A}^k$ and $\mathcal{A}^{p,q}$ as the case may be.
Moreover,
\begin{itemize}
\item
$d,\partial,\overline{\partial}$ are $\C$-linear;
\item
$d\circ d=0$, $\partial\circ\partial=0$, and $\overline{\partial}\circ\overline{\partial}=0$; and,
\item
$d(\overline{\phi})=\overline{d\phi}$, $\partial(\overline{\phi})=\overline{\overline{\partial}\phi}$, and $\overline{\partial}(\overline{\phi})=\overline{\partial\phi}$.
\end{itemize}

We say that $\omega\in\mathcal{A}^k(X)$ is {\em $d$-closed} if $d\omega=0$, and {\em $d$-exact} if $\omega=d\phi$ for some $\phi\in\mathcal{A}^{k-1}(X)$.
Since $d\circ d=0$ the exact forms are closed.
The {\em De Rham cohomology groups} are the complex vector spaces
$$H_{\operatorname{DR}}^k(X):=\frac{\{d\text{-closed }k\text{-forms}\}}{\{d\text{-exact }k\text{-forms}\}}.$$
We can relate this cohomology to the classical singular cohomology (which we assume the reader is familiar with) by integration:
Given a complex form $\omega\in\mathcal{A}^k(X)$ and a {\em $k$-simplex}
$$\phi:\Delta_k:=\left\{ (t_1,\dots,t_{k+1})\in[0,1]^{k+1}\mid\sum_{i=1}^{k+1}t_i=1\right\} \longrightarrow X,$$
it makes sense to consider $\displaystyle \int_{\phi}\omega:=\int_{\Delta_k}\phi^*\omega\in\C$.
Every complex $k$-form thereby determines a homomorphism from the free abelian group generated by the $k$-simplices (i.e., the group of {\em singular $k$-chains}) to the complex numbers.
We restrict this homomorphism to the {\em singular $k$-cycles} (those chains whose boundary is zero), and denote it by
$$\int\omega:\left\{k\text{-cycles on }X\right\}\longrightarrow\C.$$
If $\omega$ is $d$-closed then $\displaystyle \int\omega$ vanishes on boundaries by Stokes' theorem, and hence $\displaystyle\int\omega$ induces a complex-valued homomorphism on the {\em singular homology} group $\displaystyle H_k(X)=\frac{k\text{-cycles in }X}{k\text{-boundaries in }X}$.
Moreover, by Stokes' theorem again, $\displaystyle\int\omega=0$ if $\omega$ is $d$-exact.
So $\displaystyle\omega\mapsto\int\omega$ induces a homomorphism
$$H^k_{\operatorname{DR}}(X)\longrightarrow \Hom_{\mathbb{Z}}\big(H_k(U),\C\big)=H_{\operatorname{sing}}^k(X,\C)$$
where $H_{\operatorname{sing}}^k(X,\C)$ is the {\em singular cohomology} group with complex coefficients.

\medskip
\smallskip

\noindent
{\bf De Rham's Theorem} (cf. Section~4.3.2 of~\cite{voisin}){\bf .}
{\em 
The above homomorphism is an isomorphism:
$H^k_{\operatorname{DR}}(X)\cong H_{\operatorname{sing}}^k(X,\C)$.
}

\bigskip

\section{Saturation and K\"ahler manifolds}
\label{saturationsection}

One obstacle to the application of model-theoretic methods to compact complex manifolds is that the structure $\A$ is not saturated; for example, every element of every sort of $\A$ is $\emptyset$-definable.
However, for some sorts this can be seen to be an accident of the language of analytic sets in which we are working:
Suppose $V$ is a complex projective algebraic variety viewed as a compact complex manifold and consider the structure
$$\mathcal{V}:=(V;P_A\mid A\subset V^n\text{ is a subvariety over }\mathbb{Q}, n\geq 0)$$
where there is a predicate for every subvariety of every cartesian power of $V$ {\em defined over the rationals}.
It is not hard to see that $\mathcal{V}$ is saturated (it is $\omega_1$-compact in a countable language).
Moreover, Chow's theorem says that every complex analytic subset of projective space is complex algebraic, and hence a subset of a cartesian power of $V$ is definable in $\A$ if and only if it is definable (with parameters) in $\mathcal{V}$.
That is, with respect to the sort $V$, the lack of saturation in $\A$ is a result of working in too large (and redundant) a language.
This property was formalised in~\cite{sat} as follows.

\begin{Definition}
\label{es}
A compact complex variety $X$ is {\em essentially saturated} if there exists a countable sublanguage of the language of $\A$, $\mathcal{L}_0$, such that every subset of a cartesian power of $X$ that is definable in $\A$ is already definable (with parameters) in the reduct of $\A$ to $\mathcal{L}_0$.
\end{Definition}

The structure induced on $X$ by such a reduct will be saturated.

It turns out that essential saturation, while motivated by internal model-theoretic considerations, has significant geometric content.
The purpose of this section is to describe this geometric significance and to show in particular that compact K\"ahler manifolds are essentially saturated.

\medskip

We will make use of Barlet's construction of the space of compact cycles of a complex variety.
For $X$ any complex variety and $n$ a natural number, a ({\em holomorphic}) {\em $n$-cycle} of $X$ is a (formal) finite linear combination $\displaystyle Z=\sum_in_iZ_i$ where the $Z_i$'s are distinct $n$-dimensional irreducible compact analytic subsets of $X$, and each $n_i$ is a positive integer called the {\em multiplicity} of the component $Z_i$.\footnote{We hope the context will ensure that holomorphic $n$-cycles will not be confused with the singular $n$-cycles of singular homology discussed in the previous section.}
By $|Z|$ we mean the underlying set or {\em support} of $Z$, namely $\displaystyle \bigcup_iZ_i$.
We denote the set of all $n$-cycles of $X$ by $\mathcal{B}_n(X)$, and the set of all cycles of $X$ by $\displaystyle \mathcal{B}(X):=\bigcup_n\mathcal{B}_n(X)$.
In~\cite{barlet} Barlet endowed $\mathcal{B}(X)$ with a natural structure of a reduced complex analytic space.
If for $s\in \mathcal{B}_n(X)$ we let $Z_s$ denote the cycle represented by $s$, then the set
$\{(s,x):s\in\mathcal{B}_n(X), x\in |Z_s|\}$
is an analytic subset of $\mathcal{B}_n(X)\times X$.
Equipped with this complex structure, $\mathcal{B}(X)$ is called the {\em Barlet space of $X$}.
When $X$ is a projective variety the Barlet space coincides with the Chow scheme.

An cycle is called {\em irreducible} if it has only one component and that component is of multiplicity $1$.
In~\cite{campanathesis} it is shown that
$$\mathcal{B}^*(X):=\{s\in\mathcal{B}(X):Z_s\text{ is irreducible}\}$$
 is a Zariski open subset of  $\mathcal{B}(X)$.
An irreducible component of $\mathcal{B}(X)$ is {\em prime} if it has nonempty intersection with $\mathcal{B}^*(X)$.
Suppose $S$ is a prime component of the Barlet space and set
$$G_S:=\{(s,x):s\in S, x\in |Z_s|\}.$$
Then $G_S$ is an irreducible analytic subset of $S\times X$ and, if $\pi:G_S\to S$ denotes the projection map, the general fibres of $\pi$ are reduced and irreducible.
We call $G_S$ the {\em graph} of (the family of cycles parametrised) by $S$.

\begin{Fact}[cf. Theorem~3.3 of~\cite{sat}]
\label{es=compactbarlet}
A compact complex variety $X$ is essentially saturated if and only if every prime component of $\mathcal{B}(X^m)$ is compact for all $m\geq 0$.
\end{Fact}

One direction of~\ref{es=compactbarlet} is straightforward:
If every prime component of $\mathcal{B}(X^m)$ is compact, then they are all sorts in $\mathcal{A}$ and their graphs are definable in $\A$.
Consider the sublanguage $\mathcal{L}_0$ of the language of $\mathcal{A}$ made up of predicates for the graphs $G_S$ as $S$ ranges over all prime components of $\mathcal{B}(X^m)$ for all $m\geq 0$.
Then $\mathcal{L}_0$ is countable because the Barlet space has countably many irreducible components (this actually follows from Lieberman's Theorem~\ref{compact} below).
Every irreducible analytic subset of $X^m$, as it forms an irreducible cycle, is a fibre of $G_S$ for some such $S$, and hence is $\mathcal{L}_0$-definable.
By quantifier elimination for $\A$, it follows that every $\mathcal{A}$-definable subset of every cartesian power of $X$ is $\mathcal{L}_0$-definable.
So $X$ is essentially saturated.
The converse makes use of Hironaka's flattening therem and the universal property of the Barlet space.\footnote{Actually, this is done in~\cite{sat} with restricted {\em Douady spaces} (the complex analytic analogue of the Hilbert scheme) rather than Barlet spaces.
However, Fujiki has shown that if the components of the Barlet space are compact then the natural map from the Douady space to the Barlet space is proper (cf. Proposition~3.4 of~\cite{fujiki78}); and hence the components of the Douady space are also compact.}

It is in determining whether a given component of the Barlet space is compact that K\"ahler geometry intervenes.
We review the fundamentals of this theory now, and suggest~\cite{voisin} for further details.

Suppose $X$ is a complex manifold.
A {\em hermitian metric} $h$ on $X$ assigns to each point $x\in X$ a positive defnite {\em hermitian form} $h_x$ on the tangent space $T_{X,x}$.
That is, $h_x:T_{X,x}\times T_{X,x}\to\C$ satisfies:
\begin{itemize}
\item[(i)]
$h_x(-,w)$ is $\C$-linear for all $w\in T_{X,x}$,
\item[(ii)]
$h_x(v,w)=\overline{h_x(w,v)}$ for all $v,w\in T_{X,x}$,
and
\item[(iii)]
$h_x(v,v)>0$ for all nonzero $v\in T_{X,x}$.
\end{itemize}
Note that $h_x$ is $\R$-bilinear and that $h_x(v,-)$ is $\C$-antilinear for all $v\in T_{X,x}$.
Moreover this assignment should be smooth:
Given a co-ordinate system $(z,U)$ on $X$, a hermitian metric is represented on $U$ by
$$h=\sum_{i,j=1}^nh_{ij}dz_i\otimes d\overline{z_j}$$
where $h_{ij}:U\to \C$ are smooth functions and $dz_i\otimes d\overline{z_j}$ is the map that takes a pair of tangent vectors $(v,w)$ to the complex number $dz_i(v)d\overline{z_j}(w)$.

A hermitian metric encodes both a riemannian and a symplectic structure on $X$.
The real part of $h$, $\re(h_x):T_{X,x}\times T_{X,x}\to\R$, is positive definite, symmetric, and $\R$-bilinear.
That is, $\re(h)$ is a riemannian metric on $X$.
On the other hand, the imaginary part, $\im(h_x):T_{X,x}\times T_{X,x}\to\R$, is an alternating $\R$-bilinear map.
So $\im(h)$ is a real $2$-form on $X$.
Moreover, if in a co-ordinate system $(z,U)$ we have $\displaystyle h=\sum_{i,j=1}^nh_{ij}dz_i\otimes d\overline{z_j}$, then a straightforward calculation shows that
$\displaystyle \im(h)=-\frac{i}{2}\sum_{i,j=1}^nh_{ij}dz_i\wedge d\overline{z_j}$.
So as a complex $2$-form on $X$, $\im(h)$ is of type $(1,1)$.

The assignment $h\mapsto -\im(h)$ is a bijection between hermitian metrics and positive real $2$-forms of type $(1,1)$ on $X$.
We call $\omega:=-\im(h)$ the {\em K\"ahler form} associated to $h$.
A hermitian metric is a {\em K\"ahler metric} if its K\"ahler form is $d$-closed.
A complex manifold is a {\em K\"ahler manifold} if it admits a K\"ahler metric.
\begin{Example}
The {\em standard} K\"ahler metric on $\C^n$ is given by $\displaystyle \sum_{i=1}^n dz_i\otimes d\overline{z_i}$.
\end{Example}

\begin{Example}
Every complex manifold admits a hermitian metric (though not necessarily a K\"ahler one).
Indeed, given any complex manifold $X$, a cover $\mathcal{U}=(z^\iota,U_\iota)_{\iota\in I}$ by co-ordinate systems, and a partition of unity $\rho=(\rho_\iota)_{\iota\in I}$ subordinate to $\mathcal{U}$, $\displaystyle h:=\sum_{\iota\in I}\rho_\iota\left(\sum_{i=1}^ndz_i^\iota\otimes d\overline{z_i^\iota}\right)$ is a hermitian metric on $X$.
\end{Example}

\begin{Example}[Fubini-Study]
\label{fs}
Let $[z_0;\dots;z_n]$ be complex homogeneous co-ordinates for $\Pp_n$.
For each $i=0,\dots,n$ let $U_i$ be the affine open set defined by $z_i\neq 0$.
Let $F_i:U_i\to\R$ be the smooth function given by $\displaystyle \operatorname{log}\left(\frac{|z_0|^2+\cdots+|z_n|^2}{|z_i|^2}\right)$.
Then $i\partial\overline{\partial}F_i\in F_X^{1,1}(U_i)$ is real-valued.
Moreover, for all $j=0,\dots,n$, $i\partial\overline{\partial}F_i$ agrees with  $i\partial\overline{\partial}F_j$ on $U_i\cap U_j$.
Hence, the locally defined forms $i\partial\overline{\partial}F_0,\dots,i\partial\overline{\partial}F_n$ patch together and determine a global complex $2$-form on $\Pp_n$.
It is real-valued, of type $(1,1)$, and $d$-closed.
The associated K\"ahler metric is called the {\em Fubini-Study} metric on $\Pp_n$.
\end{Example}

Suppose $h$ is a hermitian metric on $X$.
There is strong interaction between the K\"ahler form $\omega$ and the riemannian metric $\re(h)$.
This is encapsulated in Wirtinger's formula for the volume of compact submanifolds of $X$.
When we speak of the {\em volume} of a submanifold of $X$ with respect to $h$, denoted by $\vol_h$, we actually mean the riemannian volume with respect to $\re(h)$.

\medskip

\noindent
{\bf Wirtinger's Formula} (cf. Section~3.1 of~\cite{voisin}){\bf .}
{\em 
If $Z\subset X$ is a compact complex submanifold of dimension $k$, then
\begin{eqnarray*}
\label{wf}
\operatorname{vol}_h(Z) & = & \int_Z\omega^k
\end{eqnarray*}
where $\omega^k=\omega\wedge\cdots\wedge\omega$ is the $k$th exterior power of $\omega$.
}

\medskip

Note that $Z$ is of real-dimension $2k$ and $\omega^k$ is a real $2k$-form on $X$, and hence it makes sense to integrate $\omega^k$ along $Z$.
If $X=\mathbb{P}_n$ and $h$ is the Fubini-Study metric of Example~\ref{fs}, then for any algebraic subvariety $Z\subseteq \mathbb{P}_n$, $\vol_h(Z)$ is the degree of $Z$.

For possibly singular complex analytic subsets  $Z\subset X$ (irreducible, compact, dimension $k$), Wirtinger's formula can serve as the {\em definition} of volume; it agrees with the volume of the regular locus of $Z$.
More generally, if $\displaystyle Z=\sum_in_iZ_i$ is a $k$-cycle of $X$, then the volume of $Z$ with respect to $h$ is
$\displaystyle \vol_h(Z):=\sum_i n_i\vol_h(Z_i)$.

Taking volumes of cycles induces a function $\vol_h:\mathcal{B}(X)\to\mathbb{R}$ given by
$$\vol_h(s):=\vol_h(Z_s).$$
The link between hermitian geometry and saturation comes from the following striking fact.

\begin{Theorem}[Lieberman~\cite{lieberman}]
\label{compact}
Suppose $X$ is a compact complex manifold equipped with a hermitian metric $h$, and $W$ is a subset of $\mathcal{B}_k(X)$.
Then $W$ is relatively compact in $\mathcal{B}_k(X)$ if and only if $\vol_h$ is bounded on $W$.
\end{Theorem}

\begin{proof}[Sketch of proof]
Wirtinger's formula tells us that $\vol_h$ is computed by integrating $\omega^k$ over the fibres of a morphism.\footnote{To be more precise, given a component $S$ of $\mathcal{B}_k(X)$, one considers the {\em universal cycle} $Z_S$ on $S\times X$ whose fibre at $s\in S$ is the cycle $Z_s$.
Then integrating $\omega^k$ over the fibres of $Z_S\to S$ gives us $\vol_h$ on $S$.}
It follows that $\vol_h$ is continuous on $\mathcal{B}_k(X)$ and hence is bounded on any relatively compact subset.

The converse relies on Barlet's original method of constructing the cycle space.
We content ourselves with a sketch of the ideas involved.
First, let $K(X)$ denote the space of closed subsets of $X$ equipped with the Hausdorff metric topology.
So given closed subsets $A,B\subset X$,
$$\operatorname{dist}(A,B):=\frac{1}{2}\left[\operatorname{max}\{\operatorname{dist}_h(A,b):b\in B\} \ + \ \operatorname{max}\{\operatorname{dist}_h(a,B):a\in A\}\right].$$
Now suppose $W\subset\mathcal{B}(X)$ is a subset on which $\vol_h$ is bounded.
Given a sequence $(s_i:i\in\mathbb{N})$ of points in $W$ we need to find a convergent subsequence.
Consider the sequence $(|Z_{s_i}|:i\in\mathbb{N})$ of points in $K(X)$.
Since $X$ is compact, so is $K(X)$, and hence there exists a subsequence $(|Z_{s_i}|:i\in I)$ which converges in the Hausdorff metric topology to a closed set $A\subset X$.
Since $\vol_h(|Z_{s_i}|)$ is bounded on this sequence, a theorem of Bishop's~\cite{bishop} implies that $A$ is in fact complex analytic.
Now, by Barlet's construction, the topology on $\mathcal{B}(X)$ is closely related to the Hausdorff topology on $K(X)$.
In particular, it follows from the fact that $(|Z_{s_i}|:i\in I)$ converges in the Hasdorff topolgy to a complex analytic subset $A\subset X$,
that  some subsequence of $(s_i:i\in I)$ converges to a point $t\in\mathcal{B}(X)$ with $|Z_t|=A$.
In particular, $(s_i:i\in\mathbb{N})$ has a convergent subsequence.
Hence $W$ is relatively compact.
\end{proof}

\begin{Corollary}[Lieberman~\cite{lieberman}]
\label{kahler-compactbarlet}
If $X$ is a compact K\"ahler manifold then the prime components of $\mathcal{B}(X)$ are compact.
\end{Corollary}

\begin{proof}[Sketch of proof]
Let $h$ be a K\"ahler metric on $X$.
The $d$-closedness of the K\"ahler form $\omega=-\im(h)$ will imply by Wirtinger's formula that $\vol_h$ is constant on the components of the Barlet space.
We sketch the argument for this here, following Proposition~4.1 of Fujiki~\cite{fujiki78}.
Fix a prime component $S$ of $\mathcal{B}_k(X)$.
By continuity of $\vol_h$, we need only show that for sufficiently general points $s,t\in S$, $\vol_h(Z_s)=\vol_h(Z_t)$.
Let $G_S\subset S\times X$ be the graph of the cycles paramatrised by $S$, and let $\pi_X:G_S\to X$ and $\pi_S:G_S\to S$ be the natural projections.
We work with a prime component so that for general $s,t\in S$, the fibres of $G_S$ over $s$ and $t$ are the reduced and irreducible complex analytic subsets $Z_s$ and $Z_t$.
Now let $I$ be a piecewise real analytic curve in $S$ connecting $s$ and $t$.
For the sake of convenience, let us assume that there is only one piece: so we have a real analytic embedding $h:[0,1]\to S$ with $h(0)=s$, $h(1)=t$ and $I=h\big([0,1]\big)$.
Consider the semianalytic set $R:=\pi_S^{-1}(I)\subset G_S$.
Given the appropriate orientation we see that the boundary $\partial R$ of $R$ in $G_S$ is $\pi_S^{-1}(s)-\pi_S^{-1}(t)$.
Note that $\pi_X$ restricts to an isomorphism between $\pi_S^{-1}(s)$ and $Z_s$ (and similarly for $\pi_S^{-1}(t)$ and $Z_t$).
Also, if $\pi_X^*(\omega^k)$ is the pull-back of $\omega^k$ to $G_S$, then $d\pi_X^*(\omega^k)=0$ since $d\omega=0$.
Using a semianalytic version of Stokes' theorem (see, for example, Herrera~\cite{herrera}), we compute
\begin{eqnarray*}
0
& = &
\int_Rd\pi_X^*(\omega^k)\\
& = &
\int_{\partial R}\pi_X^*(\omega^k)\\
& =&
\int_{\pi_S^{-1}(s)}\pi_X^*(\omega^k)-\int_{\pi_S^{-1}(t)}\pi_X^*(\omega^k)\\
& = &
\int_{Z_s}\omega^k-\int_{Z_t}\omega^k\\
& = &
\vol_h(Z_s)-\vol_h(Z_t),
\end{eqnarray*}
We have shown that $\vol_h$ is constant on $S$, and hence $S$ is compact by Theorem~\ref{compact}.
\end{proof}

If $X$ is K\"ahler then so is $X^m$ for all $m>0$.
Hence, from Corollary~\ref{kahler-compactbarlet} and Fact~\ref{es=compactbarlet} we obtain:

\begin{Corollary}
\label{kahler-es}
Every compact K\"ahler manifold is essentially saturated.
\end{Corollary}

A complex variety is said to be of {\em K\"ahler-type} if it is the holomorphic image of a compact K\"ahler manifold.
The class of all complex varieties of K\"ahler-type is denoted by $\mathcal{C}$, and was introduced by Fujiki~\cite{fujiki78}.
This class is preserved under cartesian products and bimeromorphic equivalence.
Many of the results for compact K\"ahler manifolds discussed above extend to complex varieties of K\"ahler-type.
In particular, K\"ahler-type varieties are essentially saturated (see Lemma~2.5 of~\cite{sat} for how this follows from Corollary~\ref{kahler-es} above), and their Barlet spaces have compact components which are themselves again of K\"ahler-type.

Model-theoretically we can therefore view $\mathcal{C}$ as a many-sorted structure in the language where there is a predicate for each $G_S$ as $S$ ranges over all prime components of the Barlet space of each cartesian product of sorts.
We call this the {\em Barlet language}.\footnote{This is in analogy with the {\em Douady language} from Definition~4.3 of~\cite{sat}.}
Note that every analytic subset of every cartesian product of K\"ahler-type varieties is definable (with parameters) in this language, and so we are really looking at the full induced structure on $\mathcal{C}$ from $\mathcal{A}$.
Moreover, when studying the models of $\operatorname{Th}(\mathcal{C})$, we may treat $\mathcal{C}$ as a ``universal domain'' in the sense that we may restrict ourselves to definable sets and types in $\mathcal{C}$ itself.
This is for the following reason:
Fix some definable set $F$ in an elementary extension $\mathcal{C}'$.
So there will be some sort $X$ of $\cal C$ such that $F$ is a definable subset of the nonstandard $X'$.
Essential saturation implies that there is a countable sublanguage $\mathcal{L}_0$ such that $X|_{\mathcal{L}_0}$ is saturated and every definable subset of $X^n$ in $\mathcal{A}$ is already definable in $X|_{\mathcal{L}_0}$ (with parameters).
In particular, $F$ is definable in $X|_{\mathcal{L}_0}$ over some paramaters, say $b$, in $X'$.
The $\mathcal{L}_0$-type of $b$ is realised in $X$, by say $b_0$.
Let $F_{0}$ be defined in $X'|_{\mathcal{L}_{0}}$ over $b_{0}$ in the same way as $F$ is defined over $b$.
Then in $X'|_{\mathcal{L}_{0}}$ there is an automorphism taking $F$ to $F_{0}$.
So in so far as any structural properties of $F$ are concerned we may assume it is $\mathcal{L}_{0}$-definable {\em over $X$}.
But as $X|_{\mathcal{L}_{0}}$ is a saturated  elementary substructure of $X'|_{\mathcal{L}_{0}}$, the first order properties of $F$ are then witnessed by $F\cap X$.
The latter is now a definable set in $\mathcal{C}$.
The same kind of argument works also with types.

We can also work, somewhat more canonically, as follows:
In any given situation we will be interested in at most countably many K\"ahler-type varieties, $(X_i:i\in\mathbb{N})$, at once.
We then consider the smallest (countable) subcollection $\mathcal{X}$ of sorts from $\mathcal{C}$ containing the $X_i$'s and closed under taking prime components of Barlet spaces of cartesian products of sorts in $\mathcal{X}$.
We view $\mathcal{X}$ as a multi-sorted structure in the language where there is a predicate for each $G_S$ as $S$ ranges over all such prime components of the Barlet spaces.
We call this the {\em Barlet language of $(X_i:i\in\mathbb{N})$}.
Then $\mathcal{X}$ is saturated ($2^{\omega}$-saturated and of cardinality $2^{\omega}$), $\omega$-stable, and every analytic subset of every cartesian product of sorts in $\mathcal{X}$ is definable in $\mathcal{X}$.
Moreover, after possibly naming countably many constants, $\mathcal{X}$ has elimination of imaginaries (cf. Lemma~4.5 of~\cite{sat}).
When working with K\"ahler-type varieties we will in general pass to such countable reducts of $\mathcal{C}$ without saying so explicitly, and it is in this way that we treat $\mathcal{C}$ as a universal domain for $\Th(\mathcal{C})$.

\bigskip

\section{Holomorphic forms on K\"ahler manifolds}

%In this section we will mention, sometimes without proof, various structural properties of K\"ahler manifolds and associated cohomology groups.

In section~\ref{preliminaries} we defined the De Rham cohomology groups on any complex manifold $X$ by
$\displaystyle H_{\operatorname{DR}}^k(X):=\frac{\{d\text{-closed }k\text{-forms on }X\}}{\{d\text{-exact }k\text{-forms on }X\}}$.
There are also cohomology groups of forms coming from the operators $\overline{\partial}:\mathcal{A}^{p,q}(X)\to\mathcal{A}^{p,q+1}(X)$.
Since $\overline{\partial}\circ\overline{\partial}=0$, the $\overline{\partial}$-exact forms are $\overline{\partial}$-closed.
The {\em Dolbeaut cohomology groups} are  the complex vector spaces
$$H^{p,q}(X):=\frac{\{\overline{\partial}\text{-closed forms of type }(p,q)\}}{\{\overline{\partial}\text{-exact forms of type }(p,q)\}}.$$
We denoted by $h^{p,q}(X)$ the dimension of $H^{p,q}(X)$.

A fundamental result about K\"ahler manifolds is the following fact:

\bigskip

\noindent
{\bf Hodge decomposition} (cf. Section~6.1 of~\cite{voisin}){\bf .}
{\em 
If $X$ is a compact K\"ahler manifold then $H^{p,q}(X)$ is isomorphic to the subspace of $H^{p+q}_{\operatorname{DR}}(X)$ made up of those classes that are represented by $d$-closed forms of type $(p,q)$.
Moreover, under these isomorphisms,
$\displaystyle H_{\operatorname{DR}}^k(X)\cong\bigoplus_{p+q=k}H^{p,q}(X)$.
}

\bigskip

A consequence of Hodge decomposition is that complex conjugation, which takes forms of type $(p,q)$ to forms of type $(q,p)$, induces an isomorphism between $H^{p,q}(X)$ and $H^{q,p}(X)$.
In particular, $h^{p,q}(X)=h^{q,p}(X)$.
So for $X$ compact K\"ahler and $k$ odd, $\dim_{\mathbb{C}}H_{\operatorname{DR}}^k(X)$ -- which is called the $k$th {\em Betti number} of $X$ -- is always even.

For any complex manifold $X$, given an open set $U\subset X$, a form $\omega\in\mathcal{A}^{p,0}(U)$ is called a {\em holomorphic $p$-form} on $U$ if $\overline{\partial}\omega=0$.
The sheaf of holomorphic $p$-forms on $X$ is denoted $\Omega^p$.
Since there can be no non-trivial $\overline{\partial}$-exact forms of type $(p,0)$, the holomorphic $p$-forms make up the $(p,0)$th Dolbeaut cohomology group -- that is, $H^{p,0}(X)=\Omega^{p}(X)$.
In terms of local co-ordinates a holomorphic $p$-form is just a form $\displaystyle \omega=\sum_{|I|=p}f_Idz_I$ where each $f_I:U\to\mathbb{C}$ is holomorphic.
A straightforward calculation shows that holomorphic forms are also $d$-closed.

If $X$ has dimension $n$, then $\Omega^n$ is locally of rank $1$.
The corresponding complex line bundle is called the {\em canonical bundle} of $X$, denoted by $K_X$.
The triviality of $K_X$ is then equivalent to the existence of a nowhere zero global holomorphic $n$-form on $X$.
For $X$ K\"ahler, this is precisely the condition for $X$ to be a {\em Calabi-Yau} manifold. 

Holomorphic $2$-forms will play an important role for us.
Being of type $(2,0)$, a holomorphic $2$-form $\omega\in\Omega^2(X)$ determines a $\mathbb{C}$-bilinear map $\omega_x:T_{X,x}\times T_{X,x}\to \mathbb{C}$ for each $x\in X$.
Hence it induces a $\mathbb{C}$-linear map from $T_{X,x}$ to $\Hom_{\mathbb{C}}(T_{X,x},\mathbb{C})=\Omega^1_x$.
To say that $\omega$ is {\em non-degenerate} at $x$ is to say that this map is an isomorphism.

\begin{Definition}
An {\em irreducible hyperk\"ahler} manifold (also called {\em irreducible symplectic}) is a compact K\"ahler manifold $X$ such that (i) $X$ is simply connected and (ii) $\Omega^2(X)$ is spanned by an everywhere non-degenerate holomorphic $2$-form.
\end{Definition}

The basic properties of irreducible hyperk\"ahler manifolds can be found in Section~1 of~\cite{huybrechts99}. Such properties include: $\dim(X)$ is even, $h^{2,0}(X) = h^{0,2}(X) = 1$, and $K_X$ is trivial. (The latter is because, if $\phi$ is a holomorphic $2$-form on $X$ witnessing the hyperk\"ahler condition, and $\dim(X) = 2r$ then $\phi^{r}$ is an everwhere nonzero holomorphic $2r$-form on $X$.)  For surfaces, condition (ii) in Definition 4.1 is equivalent to the triviality of $K_X$. The so-called {\em $K3$ surfaces} are precisely the irreducible hyperk\"ahler manifolds which have dimension $2$. $K3$ surfaces have been widely studied since their introduction by Weil. 
A considerable amount of information on them can be found in~\cite{bpv}.
Irreducible hyperk\"ahler manifolds are now widely studied as higher-dimensional generalizations of $K3$ surfaces. We will see in the next section the (conjectured) role of $K3$ surfaces and higher dimensional irreducible hyperk\"ahlers in the model theory of K\"ahler manifolds.

\bigskip

Given a complex manifold $X$, a cohomology class $[\omega]\in H^k_{\operatorname{DR}}(X)$ is called {\em integral} if under the identification
$$H^k_{\operatorname{DR}}(X)=H^k_{\operatorname{sing}}(X,\mathbb{C})=H^k_{\operatorname{sing}}(X,\mathbb{Z})\otimes\mathbb{C},$$
$[\omega]$ is contained in $H^k_{\operatorname{sing}}(X,\mathbb{Z})\otimes 1$.
Equivalently, the map $\displaystyle \gamma\mapsto\int_\gamma\omega$ on real $k$-cycles is integer-valued.
Similarly, $[\omega]$ is {\em rational} if it is contained in $H^k_{\operatorname{sing}}(X,\mathbb{Z})\otimes\mathbb{Q}$ under the above identification.\footnote{We have chosen not to go through the definitions of {\em sheaf cohomology}, but for those familiar with it, $H^{k}_{\operatorname{sing}}(X,\Z)$ coincides with $H^k(X,\mathbb{Z})$, the $k$th sheaf cohomology group of $X$ with coefficients in the constant sheaf $\Z$.
Likewise for $\mathbb{Q}$, $\R$, or $\C$ in place of $\Z$.
In fact for K\"ahler manifolds, the Dolbeaut cohomology group $H^{p,q}(X)$ coincides with $H^q(X,\Omega^p)$.
In any case, the integral classes can be described as those in $H^k(X,\mathbb{Z})$ and the rational ones as those in $H^k(X,\mathbb{Q})$.}

\begin{Definition} A {\em Hodge manifold} is a compact complex manifold which admits a hermitian metric $h$ whose associated K\"ahler form $\omega$ -- which recall is a real $2$-form of type $(1,1)$ -- is $d$-closed and $[\omega]$ is integral.
\end{Definition} 

In particular, a Hodge manifold is K\"ahler.

For example, $\mathbb{P}_n(\mathbb{C})$ is a Hodge manifold.
Indeed, if $\omega\in\mathcal{A}^{1,1}(\mathbb{P}_n)$ is the K\"ahler form associated to the Fubini-Study metric (see Example~\ref{fs}), and we view $\mathbb{P}_1$ as a real $2$-cycle in $\mathbb{P}_n$, then $\displaystyle \int_{\mathbb{P}_1}\omega=\pi$.
Since the (class of) $\mathbb{P}_1$ generates $H_2(\mathbb{P}_n)$, $[\frac{1}{\pi}\omega]\in H^2_{\operatorname{DR}}(\mathbb{P}_n)$ is integral.
It follows that every projective algebraic manifold is Hodge.
A famous theorem of Kodaira (sometimes called Kodaira's embedding theorem) says the converse:

\bigskip
\noindent
{\bf Kodaira's Embedding Theorem.}
{\em
Every Hodge manifold is (biholomorphic to) a projective algebraic manifold.
}
\bigskip

A consequence of Kodaira's theorem relevant for us is:

\begin{Corollary} Any compact K\"ahler manifold with no nonzero global holomorphic $2$-forms is projective.
\end{Corollary}
\begin{proof}[Sketch of proof]
If $0=\Omega^2(X)=H^{2,0}(X)$ then also $H^{0,2}(X)=0$.
By Hodge decomposition it follows that $H^2_{\operatorname{DR}}(X)=H^{1,1}(X)$.
Now let
$$H^{1,1}(X,\mathbb{R}):=\{[\omega]:\omega\text{ is real, }d\text{-closed, type}(1,1)\}.$$
Then the set
$$C:=\{[\omega]\in H^{1,1}(X,\mathbb{R}): \omega\text{ corresponds to a K\"ahler metric}\}$$
is {\em open} in $H^{1,1}(X,\mathbb{R})$ -- 
the argument being that a small deformation of a K\"ahler metric is K\"ahler.
As $X$ is K\"ahler, $C\neq\emptyset$.
On the other hand, as $H^{1,1}(X)=H^2_{\operatorname{DR}}(X)$, we have that
$H^{1,1}(X,\mathbb{R})=H^2_{\operatorname{sing}}(X,\mathbb{Z})\otimes\mathbb{R}$, and hence $C$ must contain an element of $H^2_{\operatorname{sing}}(X,\mathbb{Z})\otimes\mathbb{Q}$.
Taking a suitable integral multiple, we obtain an integral class in $C$.
Thus $X$ is a Hodge manifold, and so projective by Kodaira's embedding theorem.
\end{proof}

In fact we will require rather a {\em relative} version, proved in a similar fashion, and attributed in~\cite{campana05} to Claire Voisin:

\begin{Corollary}
\label{voisin}
Suppose that  $f:X\to S$ is a fibration in $\mathcal{C}$,  and the generic fibre of $f$ is not projective, then there exists a global holomorphic $2$-form $\omega\in\Omega^2(X)$ whose restriction $\omega_a$ to a generic fibre $X_a$ is a nonzero global holomorphic $2$-form on $X_a$.
\end{Corollary}

\bigskip

\section{Stability theory and K\"ahler manifolds}
\label{nmd}

In this section we will discuss some outstanding problems concerning the model theory (or rather stability theory) of $\Th(\cal C)$.  One concerns identifying (up to nonorthogonality, or even some finer equivalence relation), the trivial $U$-rank $1$ types. The second is the conjecture that $\Th(\cal C)$ is {\em nonmultidimensional}.  As we shall see the problems are closely related.

Because of the results in section 3, and as discussed at the end of that section, we may treat $\mathcal{C}$ as a universal domain for $\Th(\mathcal{C})$.
The main use of this is the existence of {\em generic} points: given countably many parameters $A$ from $\mathcal{C}$, and a K\"ahler-type variety $X$, there exist points in $X$ that are not contained in any proper analytic subsets of $X$ {\em defined over $A$} in the Barlet language of $X$.
The fact that we need not pass to elementary extensions in order to find such generic points makes the model-theoretic study of $\Th(\mathcal{C})$ much more accessible than that of $\Th(\mathcal{A})$.

In \cite{pillay-thisvolume} strongly minimal sets were discussed as ``building blocks" for structures of finite Morley rank. In fact one needs a slightly more general notion, that of a {\em stationary type of $U$-rank $1$}, sometimes also called a {\em minimal type}. Let us assume for now that $T$ is a stable theory, and we work in a saturated model ${\bar M}$. A complete (nonalgebraic) type $p(x)\in S(A)$ is {\em minimal} or {\em stationary of $U$-rank $1$} if for any $B\supseteq A$, $p$ has a unique extension to a nonalgebraic complete type over $B$.  Equivalently any (relatively) definable subset of the set of realizations of $p$ is finite or cofinite.  If $X$ is a strongly minimal set defined over $A$, and $p(x)\in S(A)$ is the ``generic" type of $X$, then $p(x)$ is minimal. However not every minimal type comes from a strongly minimal set. For example take $T$ to be the theory with infinitely many disjoint infinite unary predicates $P_{i}$ (and nothing else), and take $p$ to be the complete type over $\emptyset$ axiomatized by $\{\neg P_{i}(x):i<\omega\}$. We discussed the notion of {\em modularity} of a definable set $X$ in \cite{pillay-thisvolume}. The original definition was that  $X$ is modular if for all tuples $a$, $b$ of elements of $X$, $a$ is independent from $b$ over $\acl(a)\cap \acl(b)$ (together with a fixed set of parameters over which $X$ is defined), where $\acl$ is computed in ${\bar M}^{\operatorname{eq}}$. The same definition makes sense with a type-definable set (such as the set of realizations of a complete type) in place of $X$.  So we obtain in particular the notion of a modular minimal type. The minimal type $p(x)\in S(A)$ is said to be {\em trivial} if whenever $a,b_{1},..,b_{n}$ are realizations of $p$ and $a\in \acl(A,b_{1},..,b_{n})$ then $a\in\acl(A,b_{i})$ for some $i$.  Triviality implies modularity (for minimal types).  On the other hand, if $p$ is a modular nontrivial type then $p(x)$ is nonorthogonal (see below or \cite{pillay-thisvolume}) to a minimal type $q$ which is the generic type of a definable group $G$. Assuming the ambient theory to be totally transcendental, $G$ will be strongly minimal, so $p$ will also come from a strongly minimal set.  Under the same assumption (ambient theory is totally transcendental), any nonmodular minimal type will come from a strongly minimal set. So a divergence between strongly minimal sets and minimal types is only possible for trivial types.

Let us apply these notions to $\Th({\cal C})$.  Those compact complex varieties in ${\cal C}$ whose generic type is minimal are precisely the so-called {\em simple} complex varieties.
The formal definition is that a compact complex variety $X$ in $\mathcal{C}$ is {\em simple} if it is irreducible and if $a$ is a generic point of $X$ (over some set of definition) then there is no analytic subvariety $Y$ of $X$ containing $a$ with 
$0 < \dim(Y) < \dim(X)$.  (There is  an appropriate definition not mentioning generic points, and hence also applicable to all compact complex varieties.)
We are allowing the possibility that $\dim(X) = 1$, although sometimes this case is formally excluded in the definition of simplicity.
In fact, all compact complex curves are simple.
Moreover, a projective algebraic variety is simple if and only if it is of dimension 1.
If $X$ is simple we may sometimes say ``$X$ is modular, trivial, etc." if its generic type has that property.

\begin{Example}
\label{torus}
Given a $2n$-dimensional lattice $\Lambda\leq\mathbb{C}^n$, the quotient $T=\mathbb{C}^n/\Lambda$ inherits the structure of an $n$-dimensional compact K\"ahler manifold.
Such manifolds are called {\em complex tori}.
The additive group structure on $\mathbb{C}^n$ induces a compact complex Lie group structure on $T$.
If the lattice is chosen sufficiently generally -- namely the real and imaginary parts of a $\mathbb{Z}$-basis for $\Lambda$ form an algebraically independent set over $\mathbb{Q}$ -- then it is a fact that $T$ has no proper infinite complex analytic subsets, and hence is strongly minimal.
\end{Example}

A complex torus which is {\em algebraic} (bimeromorphic with an algebraic variety) is a certain kind of complex algebraic group: an abelian variety.
So the only strongly minimal algebraic complex tori are the elliptic curves, that is the $1$-dimensional abelian varieties.

%For complex tori there may be some ambiguity in use of the expression {\em simple}.  The definition of simplicity above in the case of a complex torus $T$ says precisely that $T$ has no proper nontrivial analytic subvarieties. Another common meaning is that $T$ has no proper nontrivial analytic subgroup. For abelian varieties these two notions do not coincide. But they {\em do} for complex tori which are not algebraic.

If $p$ and $q$ are the generic types of $X$ and $Y$ respectively, then $p$ is nonorthogonal to $q$ (we might say $X$ is nonorthogonal to $Y$) if and only if there is a proper analytic subvariety $Z\subset X\times Y$ projecting onto both $X$ and $Y$.
Note that if $p$ and $q$ are minimal then $Z$ must be a {\em correspondence}: the projections $Z\to X$ and $Z\to Y$ are both generically finite-to-one.

% and such that for generic $a\in X$ there are only finitely many $b\in Y$ such that $(x,y)\in Z$ and dually.

\begin{Fact}
\label{trichotomy}
Let $p(x)$ be a minimal type over ${\cal C}$.
Let $X$ be the compact complex variety whose generic type is $p$. Then either:

(i)
$p$ is nonmodular in which case $X$ is an algebraic curve,

(ii)
$p$ is modular, nontrivial, in which case $X$ is nonorthogonal to (i.e. in correspondence with) a strongly minimal complex nonalgebraic torus (necessarily of dimension $>1$), or

(iii)
$p$ is trivial, and $\dim(X) >1$.
\end{Fact}

\begin{proof}
This is proved in \cite{pillayscanlon2000} for the more general case of ${\cal A}$.
We give a slightly different argument here.

From the truth of the {\em strong conjecture} (cf. \cite{pillay-thisvolume}) for ${\cal A}$ one deduces that if $p$ is nonmodular then $X$ is nonorthogonal to a simple algebraic variety $Y$.
$Y$ has to be of dimension  $1$.
Simplicity implies that $X$ is also of dimension $1$, and so, by the Riemann existence theorem, an algebraic curve. 

If $p$ is modular and nontrivial, then as remarked above, up to nonorthogonality $p$ is the generic type of a strongly minimal (modular) group $G$.
It is proved in \cite{pillayscanlon2000} that any such group is definably isomorphic to a (strongly minimal) complex torus $T$.
If $T$ had dimension $1$ then by the Riemann existence theorem it would be algebraic, so not modular.

Likewise in the trivial case, $X$ could not be an algebraic curve so has dimension $>1$. 
\end{proof}

So the classification or description of simple trivial compact complex varieties in ${\cal C}$ remains.
Various model-theoretic conjectures have been made in earlier papers: for example that they are strongly minimal, or even that they must be $\omega$-categorical when equipped with their canonical Barlet language (see section 3).

To understand the simple trivial compact {\em surfaces} we look to
the classification of compact complex surfaces carried out by Kodaira in a series of papers in the 1960's, extending the Enriques classification of algebraic surfaces.
An account of Kodaira's work appears in \cite{bpv}. In particular Table 10 in Chapter VI there is rather useful. From it we can deduce:

\begin{Proposition}
\label{trivial2=k3}
Let $X$ be a simple trivial compact complex variety of dimension $2$ which is in the class ${\cal C}$.
Then $X$ is bimeromorphic to a $K3$ surface.
Expressed otherwise, a stationary trivial minimal type in ${\cal C}$ of dimension $2$ is, up to interdefinability, the generic type of a $K3$ surface.
\end{Proposition}

\begin{proof}  The classification of Kodaira gives a certain finite collection of (abstractly defined) classes, such that every compact surface has a ``minimal model" in exactly one of the classes, in particular is bimeromorphic to something in one of the classes.
Suppose $X$ is a simple trivial surface in ${\cal C}$. Then $X$ has algebraic dimension $0$ (namely $X$ does not map holomorphically onto any algebraic variety of dimension $>0$), $X$ is not a complex torus, and $X$ has first Betti number even.  Moreover these properties also hold of any $Y$ bimeromorphic to $X$.  By looking at Table 10, Chapter VI of \cite{bpv}, the only possibility for a minimal model of $X$ is to be a $K3$ surface.
\end{proof}

Among $K3$ surfaces are (i) smooth surfaces of degree $4$ in $\Pp_{3}$, and (ii) {\em Kummer surfaces}. A Kummer surface is something obtained from a $2$-dimensional complex torus by first quotienting by the map $x\to -x$ and then taking a minimal resolution. See Chapter VIII of \cite{bpv} for more details. In particular there are algebraic $K3$ surfaces, and there are simple $K3$ surfaces which are not trivial.  However  there do exist $K3$ surfaces of algebraic dimension $0$ (that is, which do not map onto any algebraic variety) and which are not Kummer, and these {\em will be} simple and trivial (see \cite{pillay2000}). On the other hand all $K3$ surfaces are {\em diffeomorphic} (that is, isomorphic as real differentiable manifolds), and in fact they were first defined by Weil precisely as compact complex analytic surfaces diffeomorphic to a smooth quartic surface in $\Pp_{3}$.

It is conceivable, and consistent with the examples,  that the natural analogue of Proposition 5.1 holds for higher dimensions:

\bigskip

\noindent
{\bf Conjecture I.} Any simple trivial compact complex variety in ${\cal C}$ which is bimeromorphic to (or at least in correspondence with) an irreducible hyperk\"ahler manifold.
Equivalently any trivial minimal type in ${\cal C}$ is nonorthogonal to the generic type of some irreducible hyperk\"ahler manifold.

\bigskip

Note that any irreducible hyperk\"ahler manifold has even dimension. Also, as with the special case of $K3$ surfaces, there are (irreducible) hyperk\"ahlers of any even dimension which are algebraic (and hence not trivial).

\bigskip

Let us now pass to the stability-theoretic notion of {\em nonmultidimensionality}.  We start with an arbitrary complete (possibly many-sorted) stable theory $T$, and work in a saturated model ${\bar M}$ of $T$. Let $p(x)\in S(A)$, $q(y)\in S(B)$ be stationary types (over small subsets $A,B$ of ${\bar M}$). Then $p$ is said to be {\em nonorthogonal} to $q$ if there is $C\supseteq A\cup B$ and realizations $a$ of $p$, and $b$ of $q$, such that:
(i) $a$ is independent from $C$ over $A$, and $b$ is independent from $C$ over $B$, and (ii) $a$ forks with $b$ over $C$.

A stationary type $p(x)\in S(A)$ is said to be nonorthogonal to a {\em set} of parameters $B$ if $p$ is nonorthogonal to some complete type over $\operatorname{acl}(B)$.
The theory $T$ is said to be {\em nonmultidimensional} if every stationary nonalgebraic type $p(x)\in S(A)$, is nonorthogonal to $\emptyset$.
An equivalent characterization is:
\begin{itemize}
\item[$(\star)$]
Whenever $p(x,a)$ is a stationary nonalgebraic type (with domain enumerated by the possibly infinite tuple $a$), and $stp(a') = stp(a)$, then $p(x,a)$ is nonorthogonal to $p(x,a')$.
\end{itemize}

\begin{Remark} If $T$ happens to be superstable, then it suffices that $(\star)$ holds for $p(x,a)$ regular, and moreover we may assume that $a$ is a finite tuple.  If moreover $T$ has finite rank (meaning every finitary type has finite $U$-rank), then it suffices for $(\star)$ to hold for types $p(x,a)$ of $U$-rank $1$. 
\end{Remark}

A stronger condition than nonmultidimensionality is {\em unidimensionality} which says that any two stationary nonalgebraic types are nonorthogonal. This is equivalent to $T$ having exactly one model of cardinality $\kappa$ for all $\kappa > |T|$. Nonmultidimensionality was also introduced by Shelah~\cite{shelah} in connection with classifying and counting models. For totally transcendental $T$ (namely every formula has ordinal valued Morley rank), $T$ is nonmultidimensional if and only if there is some fixed cardinal $\mu_{0}$ (which will be at most $|T|$) such that {\em essentially} the models of $T$ are naturally in one-one correspondence with sequences $(\kappa_{\alpha}:\alpha <\mu_{0})$ of cardinals. When $\mu_{0}$ is finite, $T$ is called {\em finite-dimensional}. Alternatively (for $T$ superstable of finite rank) this means that there are only finitely many stationary $U$-rank $1$ types up to nonorthogonality. 

\begin{Remark}[cf.~\cite{pillaypong}]  Suppose that $T$ is superstable of finite rank and nonmultidimensional. Suppose moreover that $T$ is one-sorted and that every stationary type of $U$-rank $1$ is nonorthogonal to a type of Morley rank $1$. Then $T$ is finite-dimensional.
\end{Remark}

The following conjecture was formulated (by Thomas Scanlon and the second author) around 2000-2001. They also pointed out (in \cite{pillayscanlon2001}) that it fails for $\operatorname{Th}({\cal A})$.

\bigskip
\noindent
{\bf Conjecture II.} $\operatorname{Th}({\cal C})$ is nonmultidimensional.
\bigskip

Let $p(x,a)$ be a stationary type in ${\cal C}$ realized by $b$ say (where $a$ is a finite tuple). Then $stp(a,b)$ is the generic type of a compact complex variety $X$, $stp(a)$ is the generic type of a compact complex variety $S$ and the map $(x,y)\to x$ gives a dominant meromorphic map $f$ from $X$ to $S$, and $tp(a/b)$ is the generic type of the irreducible fibre $X_{a}$.
Without changing $p(x,a)$ we may assume that $X$ and $S$ are manifolds and that $f$ is a holomorphic submersion (so that the generic fibre of $f$ is a manifold also).
The requirement that for another realization $a'$ of $stp(a)$, $p(x,a)$ and $p(x,a')$ are nonorthogonal, becomes: for $a'$ another generic point of $S$, there is some proper analytic subset $Z$ of $X_{a}\times X_{a'}$ which projects onto both $X_{a}$ and $X_{a'}$.  So by $(\star)$ above, we obtain the following reasonably geometric account or interpretation of the nonmultidimensionality of $\Th(\cal C)$:  for any fibration $f:X\to S$ in ${\cal C}$, any two generic fibres have the feature that there is a proper analytic subset of their product, projecting onto each factor. By Remark 5.4, we may restrict to the case where the generic fibre $X_{a}$ is {\em simple}.  So we obtain:
\begin{Remark} Conjecture II is equivalent to:
{\em Whenever $f:X\to S$ is a fibration in ${\cal C}$ with generic fibre a simple compact complex manifold, then $f$ is {\em weakly isotrivial} in the sense that for any generic fibres $X_{s}$, $X_{s'}$, there is a {\em correspondence} between $X_{s}$ and $X_{s'}$.}
\end{Remark}

Of course there are other stronger conditions than weak isotriviality which a fibration $f:X\to S$ may satisfy, for example that any two generic fibres are bimeromorphic or even that any two generic fibres are biholomorphic.
If the latter is satisfied we will call the fibration {\em isotrivial}.

Let us begin a discussion of Conjecture II.
Let $f:X\to S$ be a fibration in ${\cal C}$ with simple generic fibre $X_{s}$.  By 5.2, $X_{s}$ is either (i) an algebraic curve, (ii) a simple nonalgebraic complex torus, or (iii) has trivial generic type.
In case (i), we obtain weak isotriviality (as any two algebraic curves project generically finite-to-one onto $\Pp_{1}$).  So we are reduced to cases (ii) and (iii).
Special cases of case (ii) are proved by Campana \cite{campana05}.
Assuming the truth of Conjecture I, case (iii) is also proved in \cite{campana05}.
An exposition of this work is one of the purposes of this paper and appears in the next section.

For now, we end this section with a few aditional remarks on isotriviality.

\begin{Remark}
\label{localiso}
Let $f:X\to S$ be a fibration in ${\cal C}$.
If $f$ is locally trivial in the sense that for some nonempty open subset $U$ of $S$, $X_U$ is biholomorphic to $U\times Y$ over $U$ for some compact complex variety $Y$, then $f$ is isotrivial.
\end{Remark}
\begin{proof}
By Baire category, we can find $s_{1}, s_{2}\in U$ which are mutually generic. So $X_{s_{1}}$ is isomorphic to $X_{s_{2}}$ by assumption. But $tp(s_{1},s_{2})$ is uniquely determined by the mutually genericity of $s_{1}, s_{2}$. Hence for any mutually generic $s_{1}, s_{2}\in S$, $X_{s_{1}}$ is isomorphic to $X_{s_{2}}$. Now given generic $s_{1}, s_{2}\in S$, choose $s\in S$ generic over $\{s_{1}, s_{2}\}$. So $X_{s}$ is isomorphic to each of $X_{s_{1}}$, $X_{s_{2}}$. 
\end{proof}

\begin{Remark} Suppose that $f:X\to S$ is a fibration in ${\cal C}$ whose generic fibre $X_s$ is a simple nonalgebraic complex torus.  Suppose moreover that any (some) two mutually generic fibres $X_{s}$, $X_{s'}$ are nonorthogonal. Then any two generic fibres are isomorphic (as complex tori).
\end{Remark}
\begin{proof}
Fix two mutually generic fibres $X_s$ and $X_s'$. These are both locally modular strongly minimal groups. Hence nonorthogonality implies that there is a strongly minimal subgroup $C$ of $X_{s}\times X_{s'}$ projecting onto both factors, and this induces an isogeny from $X_{s'}$ onto $X_{s}$, and thus an isomorphism (of complex tori) between $X_{s'}/A_{s'}$ and $X_{s}$ for some finite subgroup $A_{s'}$ of $X_{s'}$. Note that $A_{s'}$ is $\operatorname{acl}(s')$-definable. Now let $s_{1}, s_{2}$ be generic points of $S$. Let $s'\in S$ be generic over $\{s_{1},s_{2}\}$. So there is an isomorphism $f_{1}$ between $X_{s'}/A_{s'}$ and $X_{s_{1}}$ (for some finite, so $\operatorname{acl}(s')$-definable subgroup of $X_{s'}$) with $X_{s_{1}}$. As $s_{1}$ and $s_{2}$ have the same type over $\operatorname{acl}(s')$, we obtain an isomorphism $f_{2}$ between $X_{s'}/A_{s'}$ and $X_{s_{2}}$. Thus $X_{s_{1}}$ and $X_{s_{2}}$ are isomorphic. 
\end{proof}

\bigskip

\section{Local Torelli and the isotriviality theorem}

In this section we state and sketch the proof of a recent result of Campana~\cite{campana05}, which was motivated by and partially resolves the nonmultidimensionality conjecture for $\mathcal{C}$ discussed in the previous section.

Suppose $Y$ is a compact K\"ahler manifold and consider a {\em deformation} $f:X\to S$ -- that is, $f$ is a proper holomorphic submersion between complex manifolds $X$ and $S$ and there is a point $o\in S$ such that $X_o=Y$.
For $s$ near $o$, $X_s$ will be a compact K\"ahler manifold (see Theorem~9.23 of~\cite{voisin}).
Diffeomorphically, $f$ is locally trivial: there exists an open neighbourhood $U\subseteq S$ of $o$ such that $X_U$ is diffeomorphic to $U\times Y$ over $U$.
Letting $u$ be this diffeomorphism we have the commuting diagram:
$$\xymatrix{
U\times Y\ar[dr]\ar[rr]^u_{\text{diffeo}} && X_U\ar[dl]^{f_U}\\
& U
}$$
For each $s\in U$, the diffeomorphism $u_s:Y\to X_s$ induces an isomorphism of singular cohomology groups, and hence by De Rham's Theorem, of the De Rham cohomology groups.
In particular, we obtain a group isomorphism, $\hat{u}_s:H^2_{\operatorname{DR}}(X_s)\to H^2_{\operatorname{DR}}(Y)$.
From our discussion of De Rham's theorem in Section~\ref{preliminaries} it is not hard to see that, under the identification
$$H^2_{\operatorname{DR}}(Y)=H^2_{\operatorname{sing}}(Y,\mathbb{C})=\Hom_{\mathbb{C}}\left(H_2(Y),\mathbb{C}\right),$$
the isomorphism $\hat{u}_s:H^2_{\operatorname{DR}}(X_s)\to H^2_{\operatorname{DR}}(Y)$ is given by
$$[\omega]\longmapsto\left([\gamma]\mapsto\int_\gamma u_s^*\omega\right)$$
where $\omega$ is a $d$-closed $2$-form on $X_s$ and and $\gamma$ is a real $2$-cycle on $Y$.

Since $u$ may not be biholomorphic, $\hat{u}_s$ does not necessarily respect the Hodge decomposition of $H^2_{\operatorname{DR}}(X_s)$ and $H^2_{\operatorname{DR}}(Y)$.
Indeed, one measure of how far $u$ is from being a biholomorphic trivialisation is the {\em period map} of $Y$ for holomorphic $2$-forms (with respect to the deformation $f$):
$$p:U\to\operatorname{Grass}\left(H^2_{\operatorname{DR}}(Y)\right)$$
which assigns to each $s\in U$ the subspace $\hat{u}_s\big(H^{2,0}(X_s)\big)$.
Recall that for any complex manifold $M$, $H^{2,0}(M)$ is just the space $\Omega^2(M)$ of global holomorphic $2$-forms on $M$.

Note that if $u$ is a biholomorphism then the period map is constant on $U$ since for all $s\in U$ $\hat{u}_s\big(H^{2,0}(X_s)\big)=H^{2,0}(Y)$.

\begin{Definition}
\label{torelli}
Suppose $Y$ is a compact K\"ahler manifold.
We say that $Y$ satisfies {\em local Torelli for holomorphic $2$-forms} if the following holds:
given any deformation $f:X\to S$ of $Y$ with a local diffeomorphic trivialisation $u:U\times Y\to X_U$, if the corresponding period map is constant on $U$ then $u$ is in fact a biholomorphic trivialisation.
\end{Definition}

\begin{Example}
\label{torih}
Complex tori and irreducible hyperk\"ahler manifolds all satisfy local Torelli for holomorphic $2$-forms.
(See Theorem 5(b) of \cite{beauville} for the case of irreducible hyperk\"ahler manifolds.)
\end{Example}

We can now state the isotriviality theorem we are interested in.

\begin{Theorem}[Campana~\cite{campana05}]
\label{isotheorem}
Suppose $f:X\to S$ is a fibration where $X$ and $S$ are compact K\"ahler manifolds.
Assume that for $a\in S$ generic,
(i) $X_{a}$ is not projective,
(ii) $\dim_{\C}\Omega^2(X_a) = 1$, and
(iii) $X_{a}$ satisfies local Torelli for holomorphic $2$-forms.
Then $f$ is isotrivial.
\end{Theorem}

\begin{proof}
[Sketch of proof]
From condition (i) and Corollary~\ref{voisin} there exists a global holomorphic $2$-form $\omega\in\Omega^2(X)$ whose restriction $\omega_a$ to the generic fibre $X_a$ is a nonzero global holomorphic $2$-form on $X_a$.
Moreover  by condition (ii), $\omega_a$ spans $\Omega^2(X_a)$.

Let $U$ be an open neighbourhood of $a$ such that there is a diffeomorphic trivialisation $u:U\times X_a\to X_U$ over $U$.
We show that the corresponding period map $p:U\to\operatorname{Grass}\big(H^2_{\operatorname{DR}}(X_a)\big)$ is constant on $U$.
By local Torelli this will imply that $u$ is biholomorphic and hence, by Remark~\ref{localiso}, $f$ is isotrivial.

For any $s\in U$ let $\hat{u}_s:H^2_{\operatorname{DR}}(X_s)\to H^2_{\operatorname{DR}}(X_a)$ be the isomorphism induced by $u$ and discussed above.
We need to show that $\hat{u}_s\big(H^{2,0}(X_s)\big)=\hat{u}_t\big(H^{2,0}(X_t)\big)$ for all $s,t\in U$.
But, shrinking $U$ if necessary, $H^{2,0}(X_s)=\Omega^2(X_s)$ is spanned by the restriction $\omega_s$ of $\omega$ to $X_s$, for all $s\in U$.
Hence it suffices to show that $\hat{u}_s(\omega_s)=\hat{u}_t(\omega_t)$ for all $s,t\in U$.

Now fix a real $2$-cycle $\gamma$ on $X_a$.
Viewing $\hat{u}_s(\omega_s)$ and $\hat{u}_t(\omega_t)$ as elements of $\Hom_{\mathbb{C}}\big(H_2(X_a),\mathbb{C}\big)=H^2_{\operatorname{DR}}(X_a)$ we compute
$$
\big(\hat{u}_s(\omega_s)-\hat{u}_t(\omega_t)\big)[\gamma]
 =  \int_\gamma u_s^*\omega_s-\int_\gamma u_t^*\omega_t
 =  \int_{u_s\circ\gamma}\omega_s- \int_{u_t\circ\gamma}\omega_t
$$
Here $u_s\circ\gamma$ and $u_t\circ\gamma$ are $2$-cycles on $X_s$ and $X_t$ respectively.
Viewed as $2$-cycles on $X$ we have
$$ \int_{u_s\circ\gamma}\omega_s- \int_{u_t\circ\gamma}\omega_t = \int_{(u_s\circ\gamma-u_t\circ\gamma)}\omega.$$
But $(u_s\circ\gamma-u_t\circ\gamma)$ is the boundary of some $3$-cycle $\lambda$ on $X$.
By Stokes', $\displaystyle \int_{(u_s\circ\gamma-u_t\circ\gamma)}\omega=\int_\lambda d\omega$.
Since holomorphic forms are $d$-closed it follows that
$$\big(\hat{u}_s(\omega_s)-\hat{u}_t(\omega_t)\big)[\gamma]=0$$
for all $2$-cycles $\gamma$ on $X_a$.
That is, $\hat{u}_s(\omega_s)=\hat{u}_t(\omega_t)$ for all $s,t\in U$.
So the period map is constant on $U$ and $f$ is isotrivial.
\end{proof}

\begin{Remark}
\label{isotheorem-cases}
The hypotheses of Theorem 6.3 are valid in the following cases:
\newline
(a) The generic fibre $X_{a}$ is irreducible hyperk\"ahler and nonprojective,
\newline
(b) The generic fibre $X_{a}$ is a simple complex torus of dimension $2$.
\end{Remark}

\begin{proof}
We have already mentioned that complex tori and irreducible hyperk\"ahler satisy local Torelli for holomorphic $2$-forms.
Nonprojectivity is assumed in (a) and follows for (b) by the fact that the only simple projective varieties are curves.
Finally, $\dim_{\mathbb{C}}\Omega^2(X_a)=1$ is true of irreducible hyperk\"ahler manifolds by definition, and true of simple complex tori of dimension $2$ by the fact that dimension $2$ forces $\dim_{\mathbb{C}}\Omega^2(X_a)$ to be at most $1$ while nonprojectivity forces it to be at least $1$.
\end{proof}

Let us return to the nonmultidimensionality conjecture (Conjecture II) from section 5, bearing in mind the equivalence stated in Remark 5.6.

\begin{Corollary}
The nonmultidimensionality conjecture holds in $\Th({\cal C})$ for  surfaces.
In other words if $p(x)$ is a minimal type of dimension $1$ or $2$ over some model of $\Th({\cal C})$ then $p$ is nonorthogonal to $\emptyset$.
\end{Corollary}

\begin{proof} As discussed at the end of section 3 we may work in $\mathcal{C}$ itself.
Let $p(x) = tp(b/a)$ for $a,b$ from ${\cal C}$ and $b$ a generic point of an $a$-definable simple compact complex manifold $X_{a}$ of dimension $1$ or $2$.
We have already pointed out that in the case of dimension $1$ (i.e. of projective curves), $X_{a}$ is nonorthogonal to $X_{a'}$ whenever $stp(a) = stp(a')$.
So assume $X_a$ is a simple compact complex surface.
It is then not projective.
By Fact~\ref{trichotomy} and Proposition~\ref{trivial2=k3},
we may assume that $X_{a}$ is either a $2$-dimensional simple complex torus, or a nonprojective $K3$ surface.
So by Theorem~\ref{isotheorem} and Remark~\ref{isotheorem-cases}, $X_{a}$ is biholomorphic to $X_{a'}$ whenever $stp(a) = stp(a')$.
Hence $tp(b/a)$ is nonorthogonal to $\emptyset$.
\end{proof}

Condition (ii) of Theorem~\ref{isotheorem} seems rather strong, and indeed, Campana works with the following weaker condition:
A {\em rational Hodge substructure} of $H^2_{\operatorname{DR}}(X_a)$ is a  $\mathbb{C}$-vector subspace $V$ such that $V^{i,j}=\overline{V^{j,i}}$ where $V^{i,j}:=V\cap H^{i,j}(X_a)$, and $V=V_{\mathbb{Q}}\otimes\mathbb{C}$ where $V_{\mathbb{Q}}:=V\cap H^2_{\operatorname{sing}}(X_a,\mathbb{Q})$.
Campana says that $X_a$ is {\em irreducibile in weight $2$} if for any rational Hodge substructure $V\subseteq H^2_{\operatorname{DR}}(X_a)$, either $V^{2,0}=0$ or $V^{2,0}=H^{2,0}(X_a)$.
By a theorem of Deligne, the image of $H_{\operatorname{DR}}^2(X)$ in $H_{\operatorname{DR}}^2(X_a)$ under the restriction map is a rational Hodge substructure.
Hence, the above proof of Theorem~\ref{isotheorem} works if condition (ii) is replaced by the irreducibility of $X_a$ in weight $2$.
Campana proves that the ``general" torus of dimension $\geq 3$ is irreducible in weight $2$.
Apparently it is open whether any simple nonalgebraic torus is irreducible in weight $2$. 
This together with Conjecture I (that any simple trivial compact K\"ahler manifold is nonorthogonal to an irreducible hyperkahler manifold) are the remaining obstacles to the nonmultidimensionality conjecture for $\mathcal{C}$.

%\bibliographystyle{plain}
%\bibliography{ccs}

\end{document}